\documentclass[11pt]{article}

\usepackage{amscd,amsmath, amssymb, fancyhdr, epsfig, color, tocloft, mathtools}
\usepackage{dutchcal}

\usepackage[backref=page]{hyperref}
\renewcommand*{\backrefalt}[4]{%
	\ifcase #1 (Not cited.)%
	\or        (Cited on page~#2.)%
	\else      (Cited on pages~#2.)%
	\fi}

\hypersetup{
	colorlinks   = true,
	citecolor    = magenta,
	linkcolor    = blue,
	urlcolor     = magenta	
}

\newcommand{\version}{version 2.0,\ \ November 1, 2021}
\makeatletter
\def\x@arrow{\DOTSB\Relbar}
\def\xlongequalsignfill@{\arrowfill@\x@arrow\Relbar\x@arrow}
\providecommand{\xlongequal}[2][]{%
	\ext@arrow 0099\xlongequalsignfill@{#1}{#2}}
\def\xlongrightarrowfill@{\arrowfill@\relbar\relbar\longrightarrow}
\newcommand{\xlongrightarrow}[2][]{%
	\ext@arrow 0099\xlongrightarrowfill@{#1}{#2}}
\makeatother

\numberwithin{equation}{section}

\def\eqref#1{(\ref{#1})}

\newcommand{\g}{{\mathfrak g}}

\newcommand{\C}{{\mathbb C}}

\def\1{\sqrt{-1}\:}

\newcommand{\cntrct}                
{\hspace{2pt}\raisebox{1pt}{\text{$\lrcorner$}}\hspace{2pt}}
\newcommand{\arrow}{{\:\longrightarrow\:}}

\renewcommand{\phi}{\varphi}
\renewcommand{\epsilon}{\varepsilon}
\renewcommand{\geq}{\geqslant}

\newcommand{\Id}{\operatorname{Id}}

\newcommand{\const}{\operatorname{\text{\sf const}}}

\newcommand{\Aut}{\operatorname{Aut}}

\newcommand{\Lie}{\operatorname{Lie}}

\newcommand{\ie}{{\em id est}}

\newcounter{Mycounter}[section]
\newcounter{lemma}[section]
\renewcommand{\thelemma}{{Lemma \thesection.\arabic{lemma}}}
\newcommand{\lemma}{%
	\setcounter{lemma}{\value{Mycounter}}
	\refstepcounter{lemma}
	\stepcounter{Mycounter}
	{\noindent \bf \thelemma:\ }}

\newcounter{claim}[section]
\newcounter{sublemma}[section]
\newcounter{corollary}[section]
\renewcommand{\thecorollary}{{Corollary \thesection.\arabic{corollary}}}
\newcommand{\corollary}{%
	\setcounter{corollary}{\value{Mycounter}}
	\refstepcounter{corollary}
	\stepcounter{Mycounter}
	{\noindent \bf \thecorollary:\ }}

\newcounter{theorem}[section]
\renewcommand{\thetheorem}{{Theorem \thesection.\arabic{theorem}}}
\newcommand{\theorem}{%
	\setcounter{theorem}{\value{Mycounter}}
	\refstepcounter{theorem}
	\stepcounter{Mycounter}
	{\noindent \bf \thetheorem:\ }}

\newcounter{conjecture}[section]
\newcounter{proposition}[section]
\renewcommand{\theproposition} {{Proposition \thesection.\arabic{proposition}}}
\newcommand{\proposition}{%
	\setcounter{proposition}{\value{Mycounter}}
	\refstepcounter{proposition}
	\stepcounter{Mycounter}
	{\noindent \bf \theproposition:\ }}

\newcounter{definition}[section]
\renewcommand{\thedefinition} {{Definition~\thesection.\arabic{definition}}}
\newcommand{\definition}{%
	\setcounter{definition}{\value{Mycounter}}
	\refstepcounter{definition}
	\stepcounter{Mycounter}
	{\noindent \bf \thedefinition:\ }}

\newcounter{example}[section]
\renewcommand{\theexample}{{Example \thesection.\arabic{example}}}
\newcommand{\example}{%
	\setcounter{example}{\value{Mycounter}}
	\refstepcounter{example}
	\stepcounter{Mycounter}
	{\noindent \bf \theexample:\ }}

\newcounter{remark}[section]
\renewcommand{\theremark}{{Remark \thesection.\arabic{remark}}}
\newcommand{\remark}{%
	\setcounter{remark}{\value{Mycounter}}
	\refstepcounter{remark}
	\stepcounter{Mycounter}
	{\noindent \bf \theremark:\ }}

\newcounter{problem}[section]
\newcounter{question}[section]
\makeatletter

\@addtoreset{equation}{section}
\@addtoreset{footnote}{section}
\makeatother

\def\blacksquare{\hbox{\vrule width 5pt height 5pt depth 0pt}}
\def\endproof{\blacksquare}

\newcommand{\proof}{{\bf Proof: \ }}
\newcommand{\pstep}{{\bf Proof. Step 1: \ }}

\begin{document}
	
	\begin{center}
		{\Large\bf  Compact homogeneous locally conformally K\"ahler manifolds are Vaisman. A new proof.}\\[5mm]
		{\large
			Liviu Ornea\footnote{Liviu Ornea is  partially supported by Romanian Ministry of Education and Research, Program PN-III, Project number PN-III-P4-ID-PCE-2020-0025, Contract  30/04.02.2021},  
			Misha
			Verbitsky\footnote{Misha Verbitsky is partially supported by
				by the HSE University Basic Research Program, FAPERJ E-26/202.912/2018 
				and CNPq - Process 313608/2017-2.\\[1mm]
				
			}\\[4mm]
			
		}
		
	\end{center}

	{\small
		\hspace{0.15\linewidth}
		\begin{minipage}[t]{0.7\linewidth}
			{\bf Abstract} \\ 
An LCK manifold with potential is a complex manifold
with a K\"ahler potential on its cover, such that
any deck transformation multiplies the K\"ahler potential by 
a constant multiplier. We prove that any homogeneous
LCK manifold admits a metric with LCK potential. This is used
to give a new proof that any compact
homogeneous LCK manifold is Vaisman.
		\end{minipage}
	}

	\tableofcontents
	
	
	\section{Introduction}
	\label{_Intro_Section_}
	
A locally conformally K\"ahler (LCK)  manifold is a complex manifold $(M,I)$ equipped with a Hermitian metric $g$ (LCK metric) which is locally conformally equivalent to a K\"ahler metric.
Then the Hermitian form $\omega(x,y):=g(Ix,y)$ satisfies $d\omega=\theta \wedge \omega$, where	$\theta$ is a 1-form, called {\em the Lee form} (see \cite{do} for an introduction to the subject).

An LCK metric $g$ is called Vaisman if the Lee form is not only closed, but also {\em parallel} with respect to the Levi Civita connection of $g$. 
	
The LCK condition is conformally invariant: if $g$ is an
LCK metric on $M$, then $e^f g$ is LCK for all smooth
functions $f$ on $M$. By contrast, the Vaisman condition
is not conformally invariant. Moreover, on a compact
complex manifold, a Vaisman metric (if it exists) is
unique up to homothety in its conformal class and
coincides with the Gauduchon metric of that conformal
class. 
	
One of the most interesting problems in LCK geometry is to find sufficient conditions for an LCK manifold to admit Vaisman metrics. In this paper we prove that homogeneity is such a sufficient  condition. 
	
 An LCK manifold $(M,I,g)$ is {\em homogeneous} if it
 admits a transitive and effective action of a
 Lie group  by holomorphic isometries of the LCK
 metric. In this case, $M=G/H$, where $G$ is a connected
 Lie group, and $H$ is the
 stabilizer subgroup of $G$.

\hfill

\remark
	(i) A group $G$ as above also preserves $\omega$ and $\theta$. 
	
	(ii) Recall that the group $\Aut(M)$ of biholomorphic conformalities of a compact LCK manifold is compact. Indeed, any holomorphic conformality preserves the corresponding Gauduchon metric which is unique up to a constant. We can therefore assume that $M$ is homogeneous under the group of holomorphic conformalities.

\hfill

The following result was announced in \cite{HK1}, see also \cite{HK2}, \cite{hk3} and \cite{guan} for a subsequent discussion. A short, complete and self-contained proof appeared in \cite{gmo}. Here we  present a new proof, based on the notion of LCK structure with potential. 

\hfill

\theorem \label{hom_vai} A compact homogeneous locally conformally K\"ahler manifold is Vaisman.

\hfill

The new proof amounts in showing that the homogeneity
implies  the existence of a holomorphic
circle action whose lift to the K\"ahler  cover does not
contain only homotheties of the K\"ahler metric. This will
imply that the manifold is LCK with potential. We then
observe that the Lee form of this LCK structure with
potential has constant length, which characterizes the
Vaisman metrics among the LCK metrics with potential.

Section 2 of this paper gathers the necessary background on LCK geometry, inlcuding Vaisman and LCK with potential metrics. In section 3 we present the proof of \ref{hom_vai}. In the last section, we give new proofs for two classical results about homogeneous Vaisman manifolds.

\section{Preliminaries}

\subsection{Locally conformally K\"ahler manifolds}

\definition Let $(M,I)$ be a complex manifold, $\dim_\C M\geq 2$. It
	is called {\bf locally conformally K\"ahler} (LCK) if it
	admits a  Hermitian metric $g$ whose 
	fundamental 2-form $\omega(\cdot,\cdot):=g(\cdot, I\cdot)$
	satisfies
	\begin{equation}\label{deflck}
		d\omega=\theta\wedge\omega,\quad d\theta=0,
	\end{equation}
	for a certain closed 1-form $\theta$ called {\bf the Lee
		form}. 

\hfill

\remark 
	Definition \eqref{deflck} is equivalent to the existence
	of a  covering $\tilde M$ endowed with a K\"ahler metric $\Omega$ which is
	acted on by the deck group $\Aut_M(\tilde M)$ by homotheties. 

\hfill


\remark\label{_d_theta_} The operator $d_\theta:=d-\theta\wedge$ is called {\bf twisted de Rham operator}. It obviously satisfies $d_\theta^2=0$ and hence $(\Lambda^*M,d_\theta)$ produces a cohomology called {\bf Morse-Novikov} or {\bf twisted}. 

\hfill

The following fundamental result of I. Vaisman shows that on compact complex manifolds, K\"ahler and LCK metrics cannot coexist.

\hfill

\theorem {\bf (\cite{va_tr})}\label{vailcknotk}
	Let $(M,\omega, \theta)$ be a compact LCK manifold, not globally conformally K\"ahler (i.e. with non-exact Lee form). 
	Then $M$ does not admit a K\"ahler metric.  
\blacksquare

\subsection{Vaisman manifolds} 

\definition An LCK manifold $(M,\omega, \theta)$ is called {\bf
		Vaisman} if $\nabla\theta=0$, where $\nabla$ is the
	Levi-Civita connection of $g$. 

\hfill

\example All
diagonal Hopf manifolds are Vaisman (\cite{ov_pams}). The Vaisman compact complex
surfaces are classified in \cite{_belgun_}, see also \cite{_ovv:surf_}. 

There exist compact LCK manifolds which do not admit Vaisman metrics. Such are the LCK Inoue surfaces, \cite{_belgun_}, the Oeljeklaus-Toma manifolds, \cite{ot}, \cite{_Otiman_}, and the non-diagonal Hopf manifolds, \cite{ov_pams}, \cite{_ovv:surf_}.

\hfill

\remark\label{_canon_foli_totally_geodesic_Remark_}
	On a Vaisman manifold, the Lee field $\theta^\sharp$ and the anti-Lee field $I\theta^\sharp$ are real holomorphic ($\Lie_{\theta^\sharp}I=\Lie_{I\theta^\sharp}I=0$) and Killing ($\Lie_{\theta^\sharp}g=\Lie_{I\theta^\sharp}g=0$), see \cite{do}. Moreover $[\theta^\sharp,I\theta^\sharp]=0$. The 2-dimensional foliation $\Sigma$   they generate is called  {\bf the canonical foliation}. 
%

\hfill

According to the above Remark, the Lee field generates a flow of holomorphic isometries. The following characterization is a partial converse:

\hfill

\theorem \label{kami_or} {\bf (\cite{kor})} 
Let $(M,\omega, \theta)$ be an LCK manifold equipped with a 
holomorphic and conformal $\C$-action $\rho$ without fixed points,
which lifts to non-isometric homotheties on 
the K\"ahler cover $\tilde M$. {Then $(M,\omega, \theta)$
	is conformally equivalent with a Vaisman manifold.}
\blacksquare

\hfill

\remark
Since $\theta$ is parallel, it has constant norm and thus we can always scale the LCK metric such that $|\theta|=1.$ In this assumption, the following formula holds,  \cite{do}:
	\begin{equation}\label{dctheta}
		d\theta^c=\omega - \theta\wedge\theta^c, \quad \text{where}\quad \theta^c(X)=-\theta(IX).
	\end{equation}
	Moreover, one can see, \cite{_Verbitsky:Vanishing_LCHK_}, that the
	(1,1)-form $\omega_0:=d^c\theta$ is semi-positive
	definite, having all eigenvalues\footnote{The eigenvalues
		of a Hermitian form $\eta$  are the eigenvalues of the
		symmetric operator $L_\eta$ defined by the equation
		$\eta(x,Iy)=g(L_\eta x,y)$.} 
	positive, except one which is 0.

\subsection{LCK manifolds with potential}

\definition An LCK manifold has  {\bf LCK potential} if it
	admits a K\"ahler covering on which the K\"ahler metric
	has a global and positive  potential function $\psi$ such that
	the deck group multiplies $\psi$ by a constant. 
In this case, $M$ is called {\bf LCK manifold with
  potential}.

\hfill

\proposition\label{dc_on_pot} The LCK manifold $(M,I,g,\theta)$ is LCK with potential if and only if equation \eqref{dctheta} is satisfied.

\hfill

\definition A function $\phi\in C^\infty (M)$ is called
            {\bf $d_\theta d^c_\theta$-plurisubharmonic}
            if $\omega= d_\theta d^c_\theta(\phi)$, where
            $d_\theta$ is the twisted de Rham operator
            (\ref{_d_theta_}) and $d^c_\theta:= I d_\theta
            I^{-1}$.

\hfill



Equation \eqref{dctheta} (and hence the definition of LCK
manifolds with potential) can be translated on the LCK
manifold itself:

\hfill

\theorem {\bf (\cite[Claim
    2.8]{ov_jgp_16})}\label{_theta_pluri_iff_pluri_}
$(M,I,\theta,\omega)$ is LCK with potential if and only if
$\omega= d_\theta d^c_\theta(\psi)$ for a strictly
positive $d_\theta d^c_\theta$-plurisubharmonic function
$\psi$ on $M$.

\hfill

\remark All Vaisman manifolds are LCK manifolds with
potential: on their K\"ahler covering, the automorphic
potential is represented by the squared length of the
pull-back of the Lee form with respect to the K\"ahler
metric. 
Among the non-Vaisman examples, we mention the
non-diagonal Hopf manifolds, \cite{ov_jgp_16}. 

\hfill

We shall need the following characterizations of the Vaisman metrics among the LCK metrics with potential:

\hfill

\proposition\label{pot_gau} {\bf (\cite{ov_pams})} 
Let $(M,\omega,\theta)$ be a compact LCK manifold with
potential. Then the LCK metric is Gauduchon if and only if
$\omega_0=d^c\theta$ is semi-positive definite, and then it is Vaisman. 
Equivalently, a compact LCK manifold with potential and with
constant norm of $\theta$ is Vaisman. \blacksquare

\hfill

In the proof that we shall present, we make use of the following sufficient condition for a compact LCK manifold to admit an LCK metric with potential (compare with \ref{kami_or}).

\hfill

\theorem  \label{_S^1_potential_Theorem_}
{\bf (\cite{ov_imrn_12}, \cite[Theorem 0.4]{_Istrati:potential_})} 
Let $M$ be a compact complex manifold, equipped
with a holomorphic $S^1$-action and an LCK metric $g$
(not necessarily $S^1$-inva\-riant). Suppose that
 this $S^1$ action does not lift to an
isometric action on the K\"ahler cover of $M$. Then $g$
admits an LCK potential.

\hfill

\proof
The automorphic K\"ahler potential 
was obtained in \cite{ov_imrn_12}
by cohomological arguments, and it is not
necessarily positive. In our terminology,
``LCK potential'' is always positive.
The existence of a positive potential is implied
a posteriori using a result from \cite{_OV:Positivity_}.
N. Istrati
(\cite[Theorem 0.4]{_Istrati:potential_}) produced an
explicit form of this argument, showing
that any $S^1$-invariant LCK metric,
in assumptions of \ref{_S^1_potential_Theorem_},
admits an LCK potential.
\endproof

%
%
%

\section{Homogeneous LCK manifolds admit LCK potential} 

The idea  is to show that the homogeneity of $G/H$ 
implies  the existence of a holomorphic
circle action which is lifted to non-isometric
homotheties of the universal cover. This will
imply that the manifold is LCK with potential (\ref{_S^1_potential_Theorem_}).
We translate this potential on the manifold itself, we average it on $G$ and obtain an invariant LCK metric with potential whose Lee form has constant length, which characterizes the
Vaisman metrics among the LCK metrics with potential (\ref{pot_gau}). The manifold is thus of Vaisman type. Finally, we prove that the initial homogeneous LCK metric itself is Vaisman.

\hfill

We start with the following lemma.

\hfill

\lemma\label{_S^1_and_theta_Lemma_}
Let $(M, \omega, \theta)$ be an LCK manifold, 
and $A$ a vector field acting on $M$ by holomorphic
isometries. Assume that 
the function $\theta(A)$ is not identically zero. 
Then $A$ does not act by isometries
when lifted to a K\"ahler cover of $M$.

\hfill

\proof
Let $\tilde M\stackrel \pi \arrow M$ be a K\"ahler cover of $M$.
Then  $\theta_1:= \pi^*\theta$ 
is exact: $\theta_1= d\phi$, and the corresponding K\"ahler form
on $\tilde M$ can be written as $\tilde\omega=e^{\phi}\pi^*\omega$.
Denote by $A_1$ the lift of $A$ to $\tilde M$. Then 
$\Lie_{A_1}(\tilde\omega) = (\Lie_{A_1} \phi) \tilde\omega$,
in other words, $A_1$ acts by isometries if and only if
$\Lie_{A_1} \phi=0$. However, Cartan's formula gives 
\[ \Lie_{A_1} \phi =A_1\cntrct d\phi=d\phi(A_1)=  \theta_1(A_1).
\]
Since, by assumption, $\theta(A)\neq 0$ somewhere on $M$, we have
$\theta_1(A_1)\neq 0$ somewhere on $\tilde M$.
\endproof

\hfill

We can prove now that $M$ is of Vaisman type.
Every homogeneous manifold $M$ can be obtained 
as $M=G/H$, where $G$ acts on $M$ transitively 
by automorphisms, and $H$ is the stabilizer of the point.
In our case, $M$ is LCK, and $G$ acts on $M$
by holomorphic LCK isometries. Since
the group of isometries of a compact manifold
is compact, we may freely assume that $G$ is
compact. This will be our running assumption
from now on.

\hfill

\proposition\label{_homo_admits_Vais_Proposition_}
Let $(M, \omega, \theta)$ be a compact, homogeneous LCK
manifold, $M=G/H$, where $G$ is a compact Lie group.
Then $M$ admits a $G$-homogeneous Vaisman metric with the same
Lee form.

\hfill

\pstep Let $\g=\Lie(G)$. Choose $A\in \g$
such that $\theta(A)\neq 0$ somewhere on $M$, and let
$T\subset G$ be the closure of the
one-parametric subgroup of $G$ generated
by $e^{tA}$. Then $T$ is a compact torus,
hence $A$ can be obtained as a limit of vector
fields $A_i$ such that each one-parametric subgroup of $G$ generated
by $e^{tA_i}$ is a circle. Then for some $i$ the vector
field $A_i$ satisfies $\theta(A_1)\neq 0$ somewhere on $M$.
By \ref{_S^1_and_theta_Lemma_}, $A_i$
it cannot act by isometries on the K\"ahler cover
of $\tilde M$. By \ref{_S^1_potential_Theorem_}, $M$ admits an LCK
metric with potential.

\hfill

{\bf Step 2:}
By \ref{_theta_pluri_iff_pluri_}, the fundamental form of
an LCK  metric with potential is written as 
$\omega_\psi=d_\theta d_\theta^c(\psi)$, for a positive
function $\psi$ on $M$. Averaging on the compact group  
$G$, we obtain a metric $g_1$ with
the fundamental form
\[ \omega_1=\mathrm{Av}_G(\omega_\psi)=d_\theta d_\theta^c(\mathrm{Av}_G(\psi)).
\]
Then $g_1$ is still LCK with potential, but the potential is now
$G$-invariant, and the metric is also $G$-invariant. Since
$\theta$ is already $G$-invariant, the averaging process
does not change it. Then $\theta$ is also the Lee form of
$g_1$, and its norm  in this metric is constant:
$|\theta|_{\omega_1}=\const$. Then \ref{pot_gau} implies
that $g_1$ is Vaisman.
\endproof

%
%
%

\hfill

Finally, we can prove that the initial homogeneous LCK
metric is itself Vaisman.

\hfill

\theorem\label{_LCK_homo_Vaisman_Theorem_}
Let $(M, g, I, \omega, \theta)$ be a compact, homogeneous LCK
manifold, $M=G/H$. Then $(M, \omega, \theta)$ is Vaisman.

\hfill

\proof
By \ref{_homo_admits_Vais_Proposition_},
$M$ admits a $G$-homogeneous Vaisman metric
$\omega_1$ with the same Lee form. 
By \ref{_S^1_potential_Theorem_}, 
we can assume that $\omega$ is an LCK
metric with potential. Then \ref{_LCK_homo_Vaisman_Theorem_}
would follow if we prove that
a $G$-homogeneous metric with potential and a given
Lee form is unique up to a constant multiplier.

\hfill

\proposition
Let $M=G/H$ be a $G$-homogeneous
complex manifold, and $\omega_1, \omega_2$ two
$G$-homogeneous LCK metrics with potential
and the same Lee form. Then 
$\omega_1$ is proportional to $\omega_2$.

\hfill

\proof By \ref{_theta_pluri_iff_pluri_},
$\omega_i = d_\theta d^c_\theta(\phi_i)$.
Averaging $\phi_i$ with $G$ if necessary, 
we can assume that it is a $G$-invariant function, 
hence constant. Then 
$\omega_i = d_\theta d^c_\theta(a_i)$,
where $a_i$ are constant functions, hence
these two forms are proportional.
\endproof

\section{Homogeneous Vaisman manifolds}

For the sake of completion, we give here new proofs of some old important results concerning homogeneous Vaisman manifolds.

\hfill 

\theorem {\bf (\cite{va_torino})} \label{homo_reg} The canonical foliation of a  compact  homogeneous Vaisman manifold $(M=G/H, I, g)$ is regular.\footnote{A foliation  on a manifold  is {\em regular} if its leaf space is a manifold; and it is {\em quasi-regular} if all its leaves are compact, in which case the leaf space is an orbifold.} 

\hfill

\proof Since $M$ is compact, $\Sigma$ has at least one compact leaf (this was proven in \cite{tsu}, but follows also by results of Kato (in \cite{kato2}) or \cite{ov_pams}). Then, by homogeneity, all leaves are compact, and hence $\Sigma$  is quasi-regular. Consider the elliptic fibration $\pi:M\rightarrow M/\Sigma$. This map has at least one smooth fibre. To see this, just consider $\pi$ with values in the smooth part of $M/\Sigma$ and apply Sard-Brown's theorem. But then, again by homogeneity, all fibers are smooth. \endproof

\hfill

\remark
On homogeneous Vaisman manifolds 
the foliation $\Sigma$ is regular, but the foliation
$\langle \theta^\sharp\rangle \subset \Sigma$ generated by the
Lee field is not necessarily regular.
Indeed, consider the Hopf surface $H=\frac{(\C^2\backslash 0)}{\langle A\rangle}$,
where $A=\alpha \Id$ and $\alpha$ is a complex number such that
$\frac{\alpha}{|\alpha|}$ is not a root of unity. In this case the Lee field 
is a radial vector field on $\C^2$, and its trajectory 
applied to a 3-dimensional sphere gives a diffeomorphism
$v\arrow \frac{\alpha}{|\alpha|}$ which has infinite order.

\hfill

A consequence of \ref{homo_reg} is the following.

\hfill

\corollary \label{_compact_homo_Vaisman_b_1_Claim_}
A compact homogeneous Vaisman manifold $M$ has $b_1=1$.

\hfill

\proof
Since the canonical foliation $\Sigma$ is regular, $M$ is a $T^2$-fibration
over a homogeneous projective manifold $P$. Since any algebraic group is rational,
homogeneous projective manifolds are also rational. Using \cite[Chapter 4]{_Debarre_},
we obtain that $P$ is simply connected. 

Write the exact sequence
\[
0 \arrow H^1(P) \arrow H^1(M) \arrow
H^1(T^2)\stackrel \gamma\arrow H^2(P)
\]
and observe that the rank of $\gamma$ is 1 (see \cite{ov_imrn_10}).
Therefore, $b_1(P)=0$ implies that
$b_1(M)=1$.
\endproof

\hfill

\remark One can construct homogeneous Vaisman manifolds as in
\cite{va_torino}, starting from a compact homogeneous
projective manifold $P=G/H$ such that the action of $G$ on
$P$ can be linearized ({\ie} $P$ admits an equivariant ample
line bundle). Note that the structure of compact
homogeneous K\"ahler manifolds is clarified in
\cite{mats}: up to biholomorphisms, they are products of
flag manifolds. 

\hfill

\noindent{\bf Acknowledgment:} We are grateful to Florin Belgun for pointing out an error in the initial proof of \ref{_LCK_homo_Vaisman_Theorem_}; thanks to him, the actual proof is much simpler.


{\small

}

\hfill

{\small
	
	\noindent {\sc Liviu Ornea\\
		University of Bucharest, Faculty of Mathematics and Informatics, \\14
		Academiei str., 70109 Bucharest, Romania}, and:\\
	{\sc Institute of Mathematics ``Simion Stoilow" of the Romanian
		Academy,\\
		21, Calea Grivitei Str.
		010702-Bucharest, Romania\\
		\tt lornea@fmi.unibuc.ro,   liviu.ornea@imar.ro}
	
	\hfill

	\noindent {\sc Misha Verbitsky\\
		{\sc Instituto Nacional de Matem\'atica Pura e
			Aplicada (IMPA) \\ Estrada Dona Castorina, 110\\
			Jardim Bot\^anico, CEP 22460-320\\
			Rio de Janeiro, RJ - Brasil }\\
		also:\\
		Laboratory of Algebraic Geometry, \\
		Faculty of Mathematics, National Research University 
		Higher School of Economics,
		6 Usacheva Str. Moscow, Russia}\\
	\tt verbit@verbit.ru, verbit@impa.br }

\end{document}